\documentclass[11pt]{amsart}
\usepackage[french]{babel}
\newtheorem{theoreme}{Th\'eor\`eme}[section]
\newtheorem{lemma}[theoreme]{Lemme}
\newtheorem{corollary}[theoreme]{Corollaire}
\newtheorem{definition}[theoreme]{D\'efinition}

\newtheorem{remarque}[theoreme]{Remarque}

 \usepackage[latin1]{inputenc}
  \usepackage[T1]{fontenc}
  \usepackage{amsmath,amssymb}
  \usepackage[all]{xy}
\begin{document}
\title[Cat\'egories associ\'ees \`a la matrice $(m_{ij}=2)$]{Classification des cat\'egories associ\'ees \`a la matrice
des coefficients $(m_{ij}=2)$ d'ordre donn\'ee}
\author{Samer Allouch}
\address{Laboratoire J. A. Dieudonn\'e\\Universit\'e de Nice-Sophia Antipolis}
\thanks{Ce papier a b\'en\'efici\'e d'une aide de l'Agence Nationale de la Recherche
portant la r\'ef\'erence ANR-09-BLAN-0151-02 (HODAG). } \maketitle

\section{Introduction}
Dans le  papier \cite{Allouch3} on a \'etudi\'e l'existence d'une
cat\'egorie ayant une matrice donn\'ee, dans ce papier on va
traviller sur $M_2^n$ la matrice 2 d'ordre n dont les
coefficients sont tous \'egaux \`a 2.\\
On note par $\mathcal{C}ard(M_2^n,r)$ le cardinale  des
cat\'egories r\'eduites dans $\mathcal{C}at(M_2^n)$. Pour n=1,2 et
3 nous avons trouv\'e que
$\mathcal{C}ard(M_2^1,r)=\mathcal{C}ard(M_2^2,r)=1$ et
$\mathcal{C}ard(M_2^3,r)=5$.\\
La question abord\'ee ici est de d\'eterminer la valeur de
$\mathcal{C}ard(M_2^n,r)$ \`a $n$ donn\'ee, mais dans ce papier on
a explor\'e leurs bornes par la formule suivante:
$$
2^{[n/3]^3}/n! \leq Card(M^n_2,r) \leq 18^{C^3_n}.
$$

\section{Quelques r\'emarques sur les fl\`eches  d'une cat\'egorie de matrice $M_2^2$}

Soit $\mathcal{A}$ \label{A d'une matrice 2,2 }une cat\'egorie
associ\'ee \`a la matrice $M_2^2$, dont les objets sont not\'es
par $\{\lambda^1,\lambda^2\}$. On supposera toujours que
$\mathcal{A}$ est {\em r\'eduite}, c'est-\`a-dire qu'il n'y a pas
d'isomorphisme entre objets distincts.

On notera par $E^{i,i}:= id_{\lambda ^i}$ l'identit\'e de $\lambda
^i$, et par $F^{i,i}$ l'unique endomorphisme $F^{i,i}\in
\mathcal{A}(\lambda ^i,\lambda ^i)$ distinct de l'identit\'e
$F^{i,i}\neq E^{i,i}$.

\begin{remarque}:{\label{Rq1 matrice 2 d'ordre 2}}
S'il existe $i\in\{1,2\}$ tel que $(F^{i,i})^2=E^{i,i}$ alors on a
deux r\'{e}sultats:\begin{enumerate}
    \item
    $gF^{i,i}=g',F^{i,i}f=f'$
    \item $fg=f'g'\neq fg'=f'g$\\
\end{enumerate}
Avec $A(\lambda^i,\lambda^j) =\{g,g'\}$ et $A(\lambda^j,\lambda^i)
=\{f,f'\}$.
\end{remarque}
\textbf{En effet:}\\
\textbf{partie 1)}\\
On prend par exemple $(F^{1,1})^2=E^{1,1}$, donc on va d\'emontrer
que $gF^{1,1}=g'$ et $F^{1,1}f=f'$.\\

Si $fF^{1,1}=f'F^{1,1}=f$ par exemple alors:
\begin{eqnarray*}
g'(F^{1,1})^2&=&g'\qquad\qquad\qquad\qquad\qquad\qquad\qquad\qquad\qquad\qquad\qquad\qquad\\
&=&\big[g'F^{1,1}\big]F^{1,1}\\
&=&(g)F^{1,1}\\
&=&g.
\end{eqnarray*}
Donc $g=g'$ contradiction, alors
$gF^{1,1}\neq g'F^{1,1}$.\\
D'autre part, si $gF^{1,1}=g$ alors il y a deux cas :
\begin{enumerate}
    \item Si $fg=F^{1,1}$
    alors,
\begin{eqnarray*}
\big[fg\big]F^{1,1}&=F^{1,1}&\qquad\qquad\qquad\qquad\qquad\qquad\qquad\qquad\qquad\qquad\qquad\qquad\\
&=&f\big[gF^{1,1}\big]\\
&=&fg\\
&=&F^{1,1}.
\end{eqnarray*}
Donc $F^{1,1}=E^{1,1}$ contradiction.
\item Si $fg=F^{1,1}$
alors,
\begin{eqnarray*}
\big[fg\big]F^{1,1}&=&(F^{1,1})^2\qquad\qquad\qquad\qquad\qquad\qquad\qquad\qquad\qquad\qquad\qquad\qquad\\
&=&E^{1,1}\\
&=&f\big[gF^{1,1}\big]\\
&=&fg\\
&=&F^{1,1}.
\end{eqnarray*}
Donc $F^{1,1}=E^{1,1}$ contradiction.
\end{enumerate}
Dans les deux cas on arrive \`a une contradiction, ce qui donne:\\
$gF^{1,1}=g'$ et
La m\^{e}me id\'{e}e pour d\'emontrer que $F^{1,1}f=f'.$\\
\\
\textbf{partie 2)}\\
On a deux cas sur $fg$ suivants:
\begin{enumerate}
    \item Si $fg=F^{1,1}$
    alors;
\begin{eqnarray*}
\big[fg\big]F^{1,1}&=&(F^{1,1})^2\qquad\qquad\qquad\qquad\qquad\qquad\qquad\qquad\qquad\qquad\qquad\qquad\\
&=&f\big[gF^{1,1}\big]\\
&=&fg'.
\end{eqnarray*}
Donc $fg'=F^{1,1}\neq fg$
    \item Si $fg=F^{1,1}$
    alors;
\begin{eqnarray*}
\big[fg\big]F^{1,1}&=&(F^{1,1})^2\qquad\qquad\qquad\qquad\qquad\qquad\qquad\qquad\qquad\qquad\qquad\qquad\\
&=&E^{1,1}\\
&=&f\big[gF^{1,1}\big]\\
&=&fg'.
\end{eqnarray*}
Donc $fg'=E^{1,1}\neq fg$.
\end{enumerate}
Dans les deux cas on trouve que $f'g\neq
fg$.\\
D'autre part,  on a $fg=f'g'$, alors il y  a deux cas sur $fg$:
\begin{enumerate}
    \item Si $fg=E^{1,1}$\\
alors $g'f=gf'=F^{2,2}$
    donc;
\begin{eqnarray*}
\big[fF^{2,2}\big]g'&=&f'g'\qquad\qquad\qquad\qquad\qquad\qquad\qquad\qquad\qquad\qquad\qquad\qquad\\
&=&F^{1,1}f\big[g'\big]\\
&=&(F^{1,1})^2\\
&=&E^{1,1}.
\end{eqnarray*}
Donc $fg=f'g'$
    \item  si $gf=E^{2,2}$
    La m\^{e}me d\'{e}monstration ci-dessus.
\end{enumerate}
Dans les deux cas on arrive \`a $fg=f'g'.$

\begin{remarque}:{\label{Rq2 matrice 2 d'ordre 2}}\\
Si $(F^{i,i})^2=E^{i,i}$ alors, $\exists (f,g)\in
A(\lambda^i,\lambda^j)\times A(\lambda^j,\lambda^i)$ tel que
$fg=F^{i,i}$.\\
\end{remarque}
\textbf{En effet:}\\
Par absurde, on pose $fg=E^{i,i}$ pour toute $(f,g)\in
A(\lambda^i,\lambda^j)\times A(\lambda^j,\lambda^i)$ c'est en
contradiction avec la remarque pr\'ec\'edante voir (\ref{Rq1
matrice 2 d'ordre 2}).\\
Donc, il existe deux fl\'eches $f,g$ tel que
$fg=F^{i,i}$.\\

\begin{lemma}:\label{e=1}
$(F^{1,1})^2=E^{1,1}\Leftrightarrow (F^{2,2})^2=E^{2,2}$
\end{lemma}
\textbf{Preuve:}\\
On pose $(F^{1,1})^2=E^{1,1}$ donc on va d\'emontrer que
$(F^{2,2})^2=E^{2,2}$.\\
Par l'absurde,
soit $(F^{2,2})^2=E^{2,2}$.\\
D'apr\`es la remarque (\ref{Rq1 matrice 2 d'ordre 2}) et remarque
(\ref{Rq2 matrice 2 d'ordre 2}) alors on a:\\
Il existe $f,g$ tel que:\\
$fg=F^{1,1}$\\
$gF^{1,1}=g'$\\
$F^{1,1}f=f'$\\
$fg=f'g'=F^{1,1}$\\
$f'g=fg'=F^{1,1}$.\\
$fg=F^{1,1}\Rightarrow [gf]g=gF^{1,1}=g'$ donc $gf=F^{2,2}$ sinon $g=g'$.\\
$fg=F^{1,1}\Rightarrow f[gf]=F^{1,1}f=f'$ ce qui donne $fF^{2,2}=f'$.\\
$f'g=E^{1,1}\Rightarrow f'[gf]=f$ ce qui donne $f'F^{2,2}=f'$.\\
D'autre part,\\
$f'=fF^{2,2}=f(F^{2,2})^2=[fF^{2,2}]F^{2,2}=f'F^{2,2}=f$
contradiction.\\
Donc $(F^{2,2})^2=E^{2,2}$.\\

\begin{lemma}:
$(F^{i,i})^2=F^{i,i}$.
\end{lemma}
\textbf{Preuve:}\\
On pose, il existe i $\in\{1,2\}$ tel que $(F^{i,i})^2=E^{i,i}$.\\
Par exemple $(F^{1,1})^2=E^{1,1}$ alors d'apr\'es le lemme
pr\'ec\'edant $(F^{2,2})^2=E^{2,2}$, ce qui donne que $A$ est non
r\'eduite.\\
D'apr\`{e}s les remarques (\ref{Rq2 matrice 2 d'ordre 2}) et
(\ref{Rq1 matrice 2 d'ordre 2}) alors,  il existe deux morphismes
$f,g$ tel que
$fg=E^{1,1}$.\\
D'autre part, on pose $gf=F^{2,2}$ alors:
\begin{eqnarray*}
\big[gf\big]g&=&F^{2,2}g\qquad \qquad\qquad \qquad\qquad \qquad\qquad\qquad\qquad\qquad \qquad\qquad\qquad\qquad\\
&=&g' \textrm{\qquad\qquad voir la remarque (\ref{Rq1 matrice 2 d'ordre 2})}\\
&=&g\big[fg\big]\\
&=&gE^{1,1}\\
&=&g.
\end{eqnarray*}
Donc $g=g'$
contradiction.\\
Donc, $gf=E^{2,2}$ et $fg=E^{1,1}$, ce qui donne les deux objets
$\lambda^1$ et $\lambda^2$ sont isomorphes entre eux,
alors $\mathcal{A}$ non-r\'eduite, contradiction.\\
Donc, $(F^{i,i})^2=F^{i,i}, \forall i$.

\begin{lemma}
\label{NonE} Pour $i\neq j$ et tout couple de morphismes $f,f'\in
\mathcal{A}(\lambda ^j,\lambda ^i)$ et $g,g'\in
\mathcal{A}(\lambda ^i,\lambda ^j)$, on a:
\begin{enumerate}
    \item  $fg = F^{i,i}$
    \item  $fF^{j,j}=f'F^{j,j}$ , $F^{i,i}f=F^{i,i}f'$ ,
    $gF^{i,i}=g'F^{i,i}$ et $F^{j,j}g=F^{j,j}g'$.
\end{enumerate}

\end{lemma}
\textbf{Preuve:}\\
(1)Par absurde\\
On pose qu'il existe f,g tel que  $fg = E^{i,i}$,
comme A r\'eduite alors $gf= F^{j,j}$.\\
D'autre part,\\
$g=g(E^{i,i})=g(fg)=(gf)g=F^{j,j}g$.\\
Donc $F^{j,j}g=g$.\\
La m\^eme chose donne $fF^{j,j}=f$.\\
On a maintenant $fF^{j,j}=f$ ce qui donne $g'f= F^{j,j}$ sinon
$g'f= E^{j,j}$
\begin{eqnarray*}
g&=&E^{i,i}g\\
&=&(g'f)g \\
&=&g'(fg)\\
&=&g'E^{i,i}\\
&=&g'.
\end{eqnarray*}
Contradiction donc $g'f= F^{j,j}$.\\
Finalement, on a $fg=E^{i,i}$,$gf=g'f=F^{j,j}$,$fF^{j,j}=f$ et $F^{j,j}g=g$ .\\
Par ailleurs,\\
\begin{eqnarray*}
g&=&F^{j,j}g\\
&=&(g'f)g \\
&=&g'(fg)\\
&=&g'E^{i,i}\\
&=&g'
\end{eqnarray*}
Ce qui donne $g=g'$, contradiction.\\
Donc $\forall$ f,g, $fg=F^{i,i}$.\\
(2) on a d'apr\`es (1) $fg=f'g=F^{i,i}$ alors,\\
$fg=f'g\Rightarrow f(gf)=f'(gf)$ alors $fF^{j,j}=f'F^{j,j}$.

\begin{lemma}
\label{Fdef} Pour chaque couple $i\neq j$ il existe un unique
morphisme, qui sera not\'e $F^{i,j}\in \mathcal{A}(\lambda
^i,\lambda ^j)$ tel que $F^{i,j}=F^{j,j}g = gF^{i,i}$ pour tout
$g\in \mathcal{A}(\lambda ^i,\lambda ^j)$.
\end{lemma}
\textbf{Preuve:}\\
Par absurde,\\
Soient $g\in \mathcal{A}(\lambda ^i,\lambda ^j)$ et
$f\in\mathcal{A}(\lambda^j,\lambda^i)$, on pose $F^{j,j}g\neq
gF^{i,i}$ alors, on peut prendre par exemple: $F^{j,j}g=g'$ et
$gF^{i,i}=g$ avec $g\neq g'$.
\begin{eqnarray*}
g'&=&F^{j,j}g\\
&=&(gf)g \textrm{\qquad\qquad voir (\ref{NonE})}\\
&=&g(fg)\\
&=&gF^{i,i}\\
&=&g
\end{eqnarray*}
Ce qui donne que $g=g'$ contradiction.\\
Donc, il existe $F^{i,j}\in \mathcal{A}(\lambda ^i,\lambda ^j)$
tel que $F^{i,j}=F^{j,j}g =
gF^{i,i}$.\\
Pour $i\neq j$ on notera par $G^{i,j}$ l'unique \'el\'ement de
$\mathcal{A}(\lambda ^i,\lambda ^j)$ distinct de $F^{i,j}$.\\
Faire attention que  $F^{i,j}$ est ind\'ependant du choix de $g$
car $F^{j,j}g= F^{j,j}g'$ voir  le lemme (\ref{NonE}).

Les morphismes de $\mathcal{A}$ sont maintenant not\'es par:
\begin{eqnarray*}
\mathcal{A}(\lambda^1,\lambda^2)&=&\{ F^{1,2},G^{1,2} \}\\
\mathcal{A}(\lambda^2,\lambda^1)&=&\{  F^{2,1},G^{2,1} \}\\
\mathcal{A}(\lambda^1,\lambda^1)&=&\{  E^{1,1} = id_{\lambda ^1},F^{1,1} \}\qquad\qquad\qquad\qquad\qquad\qquad\qquad\qquad\\
\mathcal{A}(\lambda^2,\lambda^2)&=&\{  E^{2,2}= id_{\lambda
^2},F^{2,2} \}.
\end{eqnarray*}

\begin{corollary}
\label{MultiPour2objets} La table de multiplication d'une
cat\'egorie $\mathcal{A}$ r\'eduite associ\'ee \`a la matrice
$M^2_2$ est donn\'ee, avec les notations ci-dessus, par:

\begin{itemize}
    \item $E^{j,j}X^{i,j} = X^{i,j}$,\qquad $X^{i,j}E^{i,i} = X^{i,j}$
    \item $F^{j,k}X^{i,j} = F^{i,k}$,\qquad $X^{j,k}F^{i,j}=F^{i,k}$
    \item $G^{j,i}G^{i,j} = F^{i,i}$
\end{itemize}
o\`u $i,j,k\in \{ 1,2\}$ et $X$ d\'esigne une lettre $E,F,G$ parmi
les possibilit\'es suivant $i,j$ ou $k$.
\end{corollary}
\textbf{Preuve:} Ceci est une cons\'equence des lemmes
pr\'ec\'edents.
\begin{definition}
On d\'efinit $Cat(M^n_2,r,o)$ l'ensemble des classes
d'isomorphisme de cat\'egories r\'eduites ordonnees avec matrice
$M^n_2$. Donc $Card (M^n_2,r,o)=Card( Cat(M^n_2,r,o) )$.\\
Ensuite, le groupe symetrique $S_n$ avec $n!$ \'el\'ements, agit
sur cet ensemble $Cat(M^n_2,r,o)$, et l'ensemble quotient c'est
$Cat(M^n_2,r)$, qui est l'ensemble de classes d'isomorphisme de
cat\'egories r\'eduites \`a isomorphisme non-necessairement
ordonn\'e pr\`es; et $Card (M^n_2,r)=Card( Cat(M^n_2,r) )$.
\end{definition}
\begin{lemma}
$Card(M_2^2,r)=1.$
\end{lemma}
\textbf{En effet:} D'apr\`es le table de multiplication dans la
corollaire pr\'ec\'edent, on aura une seule classe des
cat\'egories
r\'eduites qui sont associes \`a $M_2^2$.\\
Donc $Card(M_2^2,r)=1.$

\section{Quelques r\'emarques sur les fl\`eches  d'une cat\'egorie de matrice $M_2^n$}

Soit maintenant $\mathcal{A}$ une cat\'egorie r\'eduite avec
ensemble d'objets $Ob(\mathcal{A})=\{ \lambda ^1,\ldots , \lambda
^n\}$ et de matrice $M_2^n$. Ceci veut dire que
$\mathcal{A}(\lambda ^i,\lambda ^j)$ a toujours $2$ \'el\'ements.

Les consid\'erations de la section pr\'ec\'edente permettent
d'\'etablir des notations pour les morphismes de $\mathcal{A}$.
D'abord on note par $E^{i,i}$ l'identit\'e de $\lambda ^i$ et par
$F^{i,i}$ l'unique morphisme distincte de $E^{i,i}$. Ensuite, pour
tout triplet $i\neq j$ on consid\`ere la sous-cat\'egorie pleine
$\mathcal{A}^{[i,j]}$ de $\mathcal{A}$ contenant les deux objets
$\lambda ^i$ et $\lambda ^j$. On note par $F^{i,j}\in
\mathcal{A}^{[i,j]}$ l'unique morphisme donn\'e par le lemme
\ref{Fdef}, et par $G^{i,j}\in \mathcal{A}^{[i,j]}$ l'unique
morphisme distinct de $F^{i,j}$.

Nous avons donc
$$
\mathcal{A}(\lambda^i,\lambda^i)=\{  E^{i,i} = id_{\lambda
^i},F^{i,i} \}
$$
et, pour $i\neq j$,
$$
\mathcal{A}(\lambda^i,\lambda^j)=\{  F^{i,j} ,G^{i,j} \} .
$$
Le corollaire \ref{MultiPour2objets} donne la table de
multiplication pour toute composition de ces morphismes qui ne
fait intervenir que deux objets.

\begin{lemma}
\label{FXF} Pour un triplet d'objets distincts $i\neq j\neq k\neq
i$, nous avons $F^{j,k}X^{i,j}=F^{i,k}$ quelque soit $X$ (i.e.
pour $X^{i,j}=F^{i,j}$ ou $G^{i,j}$), et $X^{j,k}F^{i,j}=F^{i,j}$
quelque soit $X$ (i.e. pour $X^{i,j}=F^{i,j}$ ou $G^{i,j}$).
\end{lemma}
\textbf{Preuve:}\\
Soit $X^{i,j} \in \{F^{i,j},G^{i,j}\}$ alors on a deux cas:\\
Si $X^{i,j}=F^{i,j}$ alors,
\begin{eqnarray*}
F^{j,k}X^{i,j}&=&F^{j,k}F^{i,j}\\
&=&(F^{k,k}F^{j,k})F^{i,j} \\
&=&F^{k,k}(F^{j,k}F^{i,j})\\
&=&F^{k,k}(X^{i,k})\\
&=&F^{i,k}
\end{eqnarray*}
Donc, $F^{j,k}X^{i,j}=F^{j,k}F^{i,j}=F^{i,k}$.\\
Si $X^{i,j}=G^{i,j}$ alors,
\begin{eqnarray*}
F^{j,k}X^{i,j}&=&F^{j,k}G^{i,j}\\
&=&(F^{j,k}F^{j,j})G^{i,j} \\
&=&F^{j,k}(F^{j,j}G^{i,j})\\
&=&F^{j,k}(F^{i,j})\\
&=&F^{i,k}
\end{eqnarray*}
Donc, $F^{j,k}X^{i,j}=F^{j,k}G^{i,j}=F^{i,k}$.\\
Alors, dans les 2 cas on a trouv\'e que $F^{j,k}X^{i,j}=F^{i,k},
\forall X\in\{F,G\}$.

Au vu de ce lemme, pour un triplet d'objets distincts $i\neq j\neq
k\neq i$, la multiplication
$$
\mathcal{A}(\lambda ^j,\lambda ^k)\times \mathcal{A}(\lambda
^i,\lambda ^j)\rightarrow \mathcal{A}(\lambda ^i,\lambda ^k)
$$
et enti\`erement d\'etermin\'ee par le choix entre deux cas:
$$
G^{j,k}G^{i,j}=F^{i,k} \qquad \qquad \mbox{ not\'e cas } 0, \\
G^{j,k}G^{i,j}=G^{i,k} \qquad \qquad \mbox{ not\'e cas } 1.
$$
On d\'efinit un invariant $\alpha _{\mathcal{A}}$ par:
\begin{eqnarray*}\label{fonction alfa}
 \mathcal{\alpha}:\,\,\Big\{1,...,n\}^3& \xymatrix{ \ar[rrr]  &&& } & \{ 0,1\}  \textrm{\qquad\qquad\qquad\qquad\qquad\qquad}\\
 (i,j,k)\qquad& \xymatrix{ \ar@{|->}[rrr]  &&&}&
 \alpha(i,j,k)
\end{eqnarray*}
\begin{displaymath}
\textrm{avec $\alpha(i,j,k)$=} \left\{
\begin{array}{ll}
1\qquad \textrm{si $i\neq j\neq k\neq i$ et $G^{j,k}G^{i,j}=G^{i,k}$}\\
0\qquad \textrm{sinon }.
\end{array} \right.
\end{displaymath}
Notre tache par la suite sera de d\'eterminer les conditions
n\'ecessaires et suffisantes sur une fonction $\alpha$, pour qu'il
existe une cat\'egorie $\mathcal{A}$ avec $\alpha = \alpha
_{\mathcal{A}}$.\\
\begin{lemma}\label{card(2
d'ordre 3)=5} $\mathcal{C}ard(M_2^3,r)=5$
\end{lemma}
\textbf{En effet:}\\
On remarque que la valeur de $\mathcal{C}ard(M_2^3,r)$ ne d\'epend
que la fonction $\alpha$.\\
Donc, on va compter 5 cat\'egories r\'eduites non isomorphes.\\
Soient $i,j,k\in\{1,2,3\}$ tel que $i\neq j\neq k\neq i$ alors,
\begin{enumerate}
    \item $\alpha(i,j,k)=0$ $ \forall i,j,k$  donne la premi\`ere
    cat\'egorie nom\'e $A_1$
    \item $\alpha(i,j,k)=0$ $\forall (i,j,k)\neq (1,3,2)$  donne second cat\'egorie $A_2$
    \item $\alpha(i,j,k)=0$ $\forall (i,j,k)\neq\{(1,3,2);(2,3,1)\}$ et
    $\alpha(1,3,2)=\alpha(2,3,1)=1$ donnent la troisi\'eme cat\'egorie $A_3$
    \item $\alpha(i,j,k)=0$ $\forall (i,j,k)\neq\{(1,3,2);(3,2,1)\}$ et
    $\alpha(1,3,2)=\alpha(3,2,1)=1$ donnent la quatri\'eme cat\'egorie $A_4$
    \item $\alpha(i,j,k)=0$ $\forall (i,j,k)\neq\{(1,3,2);(2,3,1);(2,1,3)\}$ et
    $\alpha(1,3,2)=\alpha(2,3,1)=(2,1,3)=1$ donnent la cinqui\'eme cat\'egorie
    $A_5$.
\end{enumerate}
\begin{definition}:\\
$\mathcal{A}$ une cat\'egorie  associ\'ee \`a $\mathcal{M}_2^n$
dont les objets sont $\{\lambda^1,...,\lambda^n\}.$\\
Soit $(i,j,k)\in\{1,...,n\}^3$ alors,\\
Si ($i\neq j\neq k$) on dit le triplet
$[\lambda^i,\lambda^j,\lambda^k]$ est un triplet distinct.\\
Si $(i= j=k)$ on dit le triplet
$[\lambda^i,\lambda^j,\lambda^k]=[\lambda^i]$ est un
triplet identit\'{e}.\\
Si $(i=j\neq k$) on dit le triplet
$[\lambda^i,\lambda^j,\lambda^k]=[\lambda^i,\lambda^k]$ est un
triplet semi-distinct.
\end{definition}
\begin{theoreme}:\label{remarque 1+1=1}\\
Si les conditions sur $\alpha$ sont v\'erifi\'ees elle correspond
a une cat\'egorie unique et toutes les cat\'egories r\'eduites
proviennent de cela. Donc la classification des cat\'egories
r\'eduites est \'{e}quivalente \`a la classification des fonctions
$\alpha$ qui satisfont aux conditions suivantes:
\begin{enumerate}
\item  soit $[i,j,k]$ un triplet distinct  alors on a
    l'expression suivante:\\
$\alpha(i,j,k)=1$ alors $\alpha(i,k,j)=0$ et $\alpha(j,i,k)=0$.
    \item Soient $i,j,k,l$ des indices distingu\'es alors on a l'\'{e}quivalence suivante:
    \begin{eqnarray*}
\alpha(i,j,k)=1&\textrm{et}&\alpha(j,l,k)=1\\
&\Updownarrow&\\
\alpha(i,j,l)=1&\textrm{et}&\alpha(i,l,k)=1.\\
\end{eqnarray*}
\end{enumerate}
\end{theoreme}

 En effet:\\
\underline{Pour (1):}\\
On a $\alpha(i,j,k)=1$, ce signifie que $G^{j,k}G^{i,j}=G^{i,k}$
alors,
\begin{eqnarray*}
G^{k,j}G^{i,k}&=&G^{k,j}(G^{j,k}G^{i,j})\qquad\qquad\qquad\qquad\qquad\qquad\qquad\qquad\qquad\qquad\qquad\qquad\\
&=&(G^{k,j}G^{j,k})G^{i,j}\\
&=&F^{j,j}G^{i,j}\\
&=&F^{i,j}
\end{eqnarray*}
Donc,$G^{k,j}G^{i,k}=F^{i,j}$ ce qui donne $\alpha(i,k,j)=0$.\\
D'autre part,
\begin{eqnarray*}
G^{i,k}G^{j,i}&=&(G^{j,k}G^{i,j})G^{j,i}\qquad\qquad\qquad\qquad\qquad\qquad\qquad\qquad\qquad\qquad\qquad\qquad\\
&=&G^{j,k}(G^{i,j}G^{j,i})\\
&=&G^{j,k}F^{j,j}\\
&=&F^{j,k}
\end{eqnarray*}
Donc,$G^{i,k}G^{j,i}=F^{j,k}$ ce qui donne $\alpha(j,i,k)=0$.\\
\underline{Pour (2):}\\
On a $\alpha(i,j,k)=1$ et $\alpha(j,l,k)=1$ alors,
$G^{j,k}G^{i,j}=G^{i,k}$ et $G^{l,k}G^{j,l}=G^{j,k}$.\\
On va d\'emontrer que $\alpha(i,j,l)=1$.\\
supposons que $\alpha(i,j,l)=0$ c.\`a.d $G^{j,l}G^{i,j}=F^{i,l}$
\begin{eqnarray*}
G^{l,k}(G^{j,l}G^{i,j})&=&G^{l,k}F^{i,l}\qquad\qquad\qquad\qquad\qquad\qquad\qquad\qquad\qquad\qquad\qquad\qquad\\
&=&F^{i,k}\\
&=&(G^{l,k}G^{j,l})G^{i,j}\\
&=&G^{j,k}G^{i,j}\\
&=&G^{i,k}
\end{eqnarray*}
Donc, $G^{i,k}=F^{i,k}$ contradiction alors, $\alpha(i,j,l)=1$.\\
La m\^eme pour d\'emontrer $\alpha(i,l,k)=1$.\\
Il reste \`a v\'erifier l'associativit\`e avec tous les valeurs de
$\alpha$.\\
Soient $[i,j,l,k]$ un quadruple distinct alors, on va d\'emontrer
que
$(G^{l,k}G^{j,l})G^{i,j}=G^{l,k}(G^{j,l}G^{i,j})$.\\
On a les 4 cas suivant:\\
\begin{tabular}{|*{3}{c|}}
    \hline
     Cas & $\alpha(i,j,l)$  & $\alpha(j,l,k)$ \\
    \hline
     1  & 1  & 0  \\
    \hline
     2  & 0  & 1 \\
    \hline
     3  & 0  & 0  \\
    \hline
     4  & 1 & 1 \\
    \hline
\end{tabular}\\
\underline{Cas1}:\\
On a $\alpha(i,j,l)=1$  et  $\alpha(j,l,k)=0$ alors
$\alpha(i,l,k)=0$ sinon alors $\alpha(i,j,l)=\alpha(i,l,k)=1$ et
d'apr\'es la condition 2 du th\'eor\`eme ci-dessus
$\alpha(j,l,k)=1$ contradiction avec l'hypoth\'ese donc
$\alpha(i,l,k)=0$.\\
D'autre part,\\
$(G^{l,k}G^{j,l})G^{i,j}=F^{j,k}G^{i,j}=F^{i,k}$\\
$G^{l,k}(G^{j,l}G^{i,j})=G^{l,k}G^{i,l}=F^{i,k}$\\
Donc, $(G^{l,k}G^{j,l})G^{i,j}=G^{l,k}(G^{j,l}G^{i,j})$.\\
\underline{Cas2} ressemble \underline{Cas1}:\\
\underline{Cas3}:\\
$(G^{l,k}G^{j,l})G^{i,j}=F^{j,k}G^{i,j}=F^{i,k}=G^{l,k}F^{i,l}=G^{l,k}(G^{j,l}G^{i,j})$.\\
\underline{Cas4}:\\
On a $\alpha(i,j,l)=1$  et  $\alpha(j,l,k)=1$ alors,
$\alpha(i,l,k)=\alpha(i,j,k)$ sinon on pose que $\alpha(i,l,k)=1$
et $\alpha(i,j,k)=0$.\\
$\alpha(i,j,l)=\alpha(i,l,k)=1$ la condition 2 donne
$\alpha(i,j,k)=1$ contradiction donc $\alpha(i,l,k)=\alpha(i,j,k)$
ce signifie que $G^{l,k}G^{i,l}=G^{j,k}G^{i,j}$.\\
Par ailleures,\\
$(G^{l,k}G^{j,l})G^{i,j}=G^{j,k}G^{i,j}=G^{l,k}G^{i,l}=G^{l,k}(G^{j,l}G^{i,j})$.\\
Soient $[i,j,l,k]$ un quadruple semi-distinct alors, on va
d\'emontrer que
$(G^{l,k}G^{j,l})G^{i,j}=G^{l,k}(G^{j,l}G^{i,j})$.\\
On a les cas suivantes:\\
Si (i=j) alors,\\
$(G^{l,k}G^{j,l})G^{i,j}=(G^{l,k}G^{i,l})F^{i,i}=X^{i,k}F^{i,i}=F^{i,k}$
avec $X\in\{F,G\}$.\\
$G^{l,k}(G^{i,l}G^{i,i})=G^{l,k}F^{i,l}=F^{i,k}$.\\
Donc, $(G^{l,k}G^{j,l})G^{i,j}=G^{l,k}(G^{j,l}G^{i,j})$.\\
Si (i=l) alors,\\
Dans ce cas on a deux cas sur $\alpha(j,i,k)$.\\
\begin{itemize}
    \item si $\alpha(j,i,k)=0$ alors $(G^{l,k}G^{j,l})G^{i,j}=(G^{i,k}G^{j,i})G^{i,j}=F^{j,k}G^{i,j}
    =F^{i,k}=G^{i,k}F^{i,i}=G^{i,k}(G^{j,i}G^{i,j})$
    \item  si $\alpha(j,i,k)=1$ alors la condition 1 du
    th\'eor\`eme ci-dessus donne que $\alpha(i,j,k)=0$, alors
$(G^{l,k}G^{j,l})G^{i,j}=(G^{i,k}G^{j,i})G^{i,j}=G^{j,k}G^{i,j}
    =F^{i,k}=G^{i,k}F^{i,i}=G^{i,k}(G^{j,i}G^{i,j})$
\end{itemize}
Donc dans ce cas on a
$(G^{l,k}G^{j,l})G^{i,j}=G^{l,k}(G^{j,l}G^{i,j})$.\\
Si (i=k) alors,\\
$(G^{l,k}G^{j,l})G^{i,j}=(G^{l,i}G^{j,l})G^{i,j}=Y^{j,i}G^{i,j}=F^{i,i}$\\
$G^{l,k}(G^{j,l}G^{i,j})=G^{l,i}(G^{j,l}G^{i,j})=G^{l,i}Z^{i,l}=F^{i,i}$\\
avec $Y?Z\in\{F,G\}$\\
Donc,$(G^{l,k}G^{j,l})G^{i,j}=G^{l,k}(G^{j,l}G^{i,j})$.\\
les qui sont rest\'ees la m$\^{e}$me id\'ees.\\
Finalement dans le cas [i,j,k,l] identit\'e bien sur il y a
l'associativit\'e.\\
\begin{lemma}
Soit $[i,j,k]$ triple semi-distinct ou identit\'e alors
$\alpha(i,j,k)=0$
\end{lemma}
\textbf{Preuve:}\\
Si i=j=k alors $G^{i,i}G^{i,i}=F^{i,i}$ c.\`a.d $\alpha(i,j,k)=0$.\\
Si i=j alors $G^{j,k}G^{i,i}=F^{i,k}$ c.\`a.d $\alpha(i,j,k)=0$.\\
Si i=k alors $G^{j,i}G^{i,j}=F^{i,i}$ c.\`a.d $\alpha(i,j,k)=0$.\\
Si j=k alors $G^{j,j}G^{i,j}=F^{i,j}$ c.\`a.d $\alpha(i,j,k)=0$.\\
\textbf{Remarque:}\\
Soit [i,j,l,k] un quadruple alors on a la formule suivante:\\
$$
\alpha (i,j,k)\alpha (j,l,k) = \alpha (i,j,l)\alpha (i,l,k).
$$
\textbf{En effet:}\\
Si $\alpha (i,j,k)=\alpha (j,l,k)=1$ alors la condition deux donne
$\alpha (i,j,l)=\alpha (i,l,k)=1$ donc,\\
$\alpha (i,j,k)\alpha (j,l,k) = \alpha (i,j,l)\alpha (i,l,k)=1$.\\
Si $\alpha (i,j,k)=\alpha (j,l,k)=0$ alors l'un de deux est
\'egale \`a 0 sinon $\alpha (i,j,l)=\alpha (i,l,k)=1$ c'est
contradiction avec la condition 2 donc,\\
$\alpha (i,j,k)\alpha (j,l,k) = \alpha (i,j,l)\alpha (i,l,k)=0$.\\
Si $\alpha (i,j,k)\neq\alpha (j,l,k)$ alors l'un de deux est
\'egale \`a 0 sinon $\alpha (i,j,l)=\alpha (i,l,k)=1$ c'est
contradiction avec la condition 2 donc,\\
$\alpha (i,j,k)\alpha (j,l,k) = \alpha (i,j,l)\alpha (i,l,k)=0$.\\
Finalement, $\alpha(i,j,k)=0$ si [i,j,k] est semi-distinct ou
identit\'e.

\section{les bornes
des cardinalit\'es} \textbf{Notation:} On va d\'efinit deux
notations:
\begin{enumerate}
    \item On veut dire par $\sigma$ la notation suivante:
\begin{displaymath}
\sigma := \lim_{n  \to  \infty} Sup\frac{\log
(Card(M^n_2,r))}{n^3}.
\end{displaymath}
    \item On note le nombre
des cat\'egories r\'eduites a isomorphisme ordonn\'ee pr\`es qui
sont associ\'ees \`{a} la matrice $M_2^n$
par $Card(M^n_2,r,o)$.\\
Par exemple: on va voir dans le lemme(5.10) $Card(M^n_2,r,o)\leq
18^{C^3_n}$; cette borne \'{e}tant atteinte dans les cas $n=2$ et
$n=3$.
\end{enumerate}
\begin{lemma}:\label{card(M,r,o)/n!,card(M,r),card(M,r,o)} on a les in\'{e}galit\'{e}s suivantes:
\begin{equation}\label{resultat}
(Card (M^n_2,r,o)/n!) \leq Card (M^n_2,r)\leq Card (M^n_2,r,o).
\end{equation}
\end{lemma}

\textbf{Preuve:} Gr\^{a}ce \`{a} chacune des d\'efinitions de
$Card (M^n_2,r)$ et de $Card (M^n_2,r,o)$, on arrive \`{a}
l'in\'equation (\ref {resultat}).
\begin{theoreme}\label{A=18}:
Le nombre de configurations de la fonction $\alpha :
(i,j,k)\mapsto 0$ ou $1$ sur ce triplet
$[\lambda^i,\lambda^j,\lambda^k]$ , est $18$.
\end{theoreme}

\textbf{Preuve:} On peut prendre ce triplet par exemple le triplet
$[\lambda^1,\lambda^2,\lambda^3]$ dans le cas de matrice $2$
d'ordre $3$.\\
On a trouv\'{e} $5$ cat\'egories non isomorphes entre elles
associ\'ees \`{a} $M_2^3$ voir le th\'eor\`{e}me(\ref{card(2
d'ordre 3)=5}).

Maintenant on a besoin de savoir combien il y a d'isomorphisme
ordonne pres, c'est a dire on ne confond pas ceux qui sont
semblables. On peut les classifier par leur classe d'isomorphisme
$A_1, A_2, A_3, A_4, A_5.$\\
$A_1$: il n'y a qu'une ici.\\

$A_2$: Pour que cela marche on doit choisir un couple $(ij)$
distinct parmi
$1,2,3$, il y a $6$ choix.\\

$A_3$: Pour que cela marche on doit choisir un couple $(ij)$ mais
$(ij)=(ji)$, vu que l'ensemble des r\'{e}ponses est $\{ (ij)$ et
$(ji)\}$ , donc il y a $3$ choix.\\

$A_4$: Pour que cela marche on doit choisir une suite de $3$
\'el\'ements
distincts, il y a $3!=6$ choix.\\

$A_5:$ c'est un cycle, qui peut aller dans un sens ou dans
l'autre, donc il y a $2$ choix:\\
$$(12) + (23) + (31)$$
$$\textrm{ou}$$
$$(13) + (32) +(21).$$
Au total on a $A= 1 + 6 + 3 + 6 + 2 = 18$ possibilit\'{e}s; le
nombre des cat\'egories non-r\'eduites a isomorphisme ordonne pres
est $18$, et c'est le nombre de fonctions $\alpha$ possible sur
$3$
indices.\\
Donc on peut dire aussi $\mathcal{C}ard (M^3_2, r,o)=18$.
\begin{lemma}:\label{borne superieur de $M_2^n$} $\mathcal{C}ard (M^n_2, r,o)\leq18^{C^3_n}.$
\end{lemma}

  En effet: Sur chaque triplet $[\lambda^i,\lambda^j,\lambda^k]$
  on a le nombre
des cat\'egories non-r\'eduites \`{a} isomorphisme ordonne pres
est $18$ voir le th\'eor\`{e}me (\ref{A=18}), alors totalemnet on
a $18^{C^3_n}$ ce qui donne $\mathcal{C}ard (M^n_2,r,o)\leq
18^{C^3_n}$.
\begin{corollary}:\label{borne superieure de card(M,r)} La borne sup\'{e}rieure de $\mathcal{C}ard
(M^n_2,r)$ est $18^{C^3_n}$ .
\end{corollary}

\textbf{Preuve:} le lemme (\ref{borne superieur de $M_2^n$}) et le
lemme (\ref{card(M,r,o)/n!,card(M,r),card(M,r,o)}) donnent:
$$\mathcal{C}ard (M^n_2,r)\leq
18^{C^3_n}$$.

 \textbf{Notation: } On consid\`{e}re $X = \{ x_1,\ldots , x_n\}$
l'ensemble (ordonn\'ee) de n objets. Soit $P_3(X)$ l'ensemble des
parties a trois elements de $X$, et $O_3(X)$ l'ensemble des
triplets distincts (avec un ordre qui peut etre different de
l'ordre de $X$). On a donc $Card(O_3(X)) = 3!Card (P_3(X))$ et
$Card (P_3(X)) = C^3_n$. Un triplet $(x_i,x_j,x_k)\in O_3(X)$ sera
not\'{e} aussi $(i,j,k)$, on a $i\neq j$, $j\neq k$ et $i\neq k$.
\begin{definition}: Soit $H\subset O_3(X)$ un sous-ensemble. On dit que $H$ est
{\em non-interferant} si: pour tout $(i,j,k)\in H$, il n'existe
pas de $(i,u,j)\in H$ ni de $(j,v,k)\in H$.
\end{definition}
\begin{lemma} Si $H\subset O_3(X)$ est un sous-ensemble non-interferant,
alors pour tout $H'\subset H$, on a $H'$ aussi
non-interf\'{e}rant.
\end{lemma}

\textbf{Preuve:} Soit $H'$ un sous ensemble de $H$, on va
d\'emontrer
$H'$ est {\em non-interf\'{e}rant}.\\
Soit $(i,j,k)\in H'\subset H$ comme $H'$ est {\em non-interferant}
alors il n'existe pas de $(i,u,j)\in H'$ ni de $(j,v,k)\in H'$,
donc $H'$ est aussi non-interf\'{e}rant.

\begin{lemma}\label{non-interferant
existe cat\'egorie}: Si $H\subset O_3(X)$ est un sous-ensemble
non-interf\'{e}rant, il existe une cat\'egorie $A_H\in M^n_2$
telle que $(i,j,k)=1$ si $(i,j,k)\in H$ et $(i,j,k)=0$ si
$(i,j,k)\not \in H$.
\end{lemma}

\textbf{ Preuve:} On construit la cat\'egorie $B$ avec les
conditions suivantes:
\begin{enumerate}
    \item  Si
$(i,j,k)\in H$ alors $\alpha(i,j,k)=1, $ et si $(i,j,k)\not\in H$
alors $\alpha(i,j,k)=0.$
    \item  Pour tout quadruplet $i,u,j,k$ alors:\\
    $\alpha(i,j,k)= 1$ et $\alpha (i,u,j)=1 \Leftrightarrow \alpha
(i,u,k)=1$ et $\alpha(u,j,k)=1$, c'est analogue a la
remarque(\ref{remarque 1+1=1})
\end{enumerate}
Alors $B$ est bien une cat\'egorie.

\begin{lemma}:
Si $H\subset O_3(X)$ est un sous-ensemble non-interf\'{e}rant,
alors il y a une cat\'egorie diff\'{e}rente (r\'eduite) $A_{H'}\in
M^n_2$ pour chaque sous-ensemble $H'\subset H$. En
cons\'{e}quence, on a
$$
\mathcal{C}ard (M^n_2,r,o)\geq 2^{Card (H)}.
$$
\end{lemma}

\textbf{En effet:} Soit $H'\subset H$ alors $H'$ est
non-interf\'{e}rant alors d'apr\`{e}s le lemme
(\ref{non-interferant existe cat\'egorie}) alors il existe une
cat\'egorie $B_{H'}$, bien sur $A_{H}\neq B_{H'}$ car $H'\subset
H$ et ce qui donne $\mathcal{C}ard (M^n_2,r,o)\geq 2^{Card (H)}.$

\begin{theoreme}: \label{borner de card(m,r,o)} On peut d\'eterminer des bornes
(voir ci-dessous) pour la cardinalit\'e de l'ensemble
$\mathcal{C}ard(M^n_2,r,o).$
\end{theoreme}

\textbf{Preuve:} Construction d'un sous-ensemble non-interferant:
supposons que $X=X_1\cup X_2\cup X_3$ est une reunion avec
$X_1\cap X_2=\emptyset$, $X_1\cap X_3=\emptyset$, et et $X_2\cap
X_3=\emptyset$(c.\`{a}.d $X$ est la reunion disjointe de $X_1$,
$X_2$, $X_3$). Alors si on pose
$$
H= \{ (i,j,k)/ x_i\in X_1, x_j\in X_2 \,\,et\,\, x_k\in X_3\}.
$$
Le sous-ensemble $H\subset O_3(X)$ est non-interferant. En effet,
les conditions impliquent d\'{e}j\'{a} que $i\neq j$, $j\neq k$ et
$i\neq k$, aussi qu'il ne peut pas y avoir d'elements $(i,u,j)$ ni
de $(j,v,k)$ dans $H$.
\\
Cet ensemble a $Card(H)= Card(X_1)Card(X_2)Card(X_3)$.
\\
On peut prendre par exemple:\\
$Card(X_1)=[n/3], Card (X_2)= [n/3]$ et $Card (X_3) = n-2[n/3]
\geq [n/3]$.\\
Donc $Card (H)\geq [n/3]^3$ qui est de l'ordre de $n^3/27$. On a
la borne inf\'{e}rieure:
$$
\mathcal{C}ard (M^n_2,r,o) \geq 2^{[n/3]^3}.
$$
Donc
$$
\log (\mathcal{C}ard (M^n_2,r,o))\geq n^3 \log (2)/27 .
$$

La borne sup\'{e}rieure est ${18}^{C^3_n}$ voir le lemme
(\ref{borne superieur de $M_2^n$})c.\`{a}.d:
$$
\mathcal{C}ard (M^n_2,r,o)\leq {18}^{C^3_n}.
$$
Et $C^3_n$ est de l'ordre de $n^3 / 6$, donc
$$
\log (\mathcal{C}ard (M^n_2,r,o))\leq n^3 \log (18) / 6 .
$$
On a donc un encadrement de $\log (\mathcal{C}ard (M^n_2,r,o))$ ou
les deux termes ont un ordre de croissance de $n^3$.
\begin{lemma}: On peut d\'{e}terminer des bornes inf\'erieures et sup\'erieures
(dont les valeurs ci-dessous) pour $\mathcal{C}ard (M^n_2,r)$.
\end{lemma}

\textbf{En effet:} Le lemme
(\ref{card(M,r,o)/n!,card(M,r),card(M,r,o)}) et le th\'eor\`{e}me
(\ref{borner de card(m,r,o)})donnet la borne sup\'{e}rieur de
$Card(M^n_2,r)$ qui est:
$$
2^{[n/3]^3}/n! \leq Card(M^n_2,r,o)/n!\leq Card(M^n_2,r).
$$
D'autre part, on a d'apr\`{e}s le corollaire (\ref{borne
superieure de
card(M,r)}) $Card(M^n_2,r)\leq 18^{C^3_n}.$\\
Donc on peut dire:
$$
2^{[n/3]^3}/n! \leq Card(M^n_2,r) \leq 18^{C^3_n}.
$$
\begin{theoreme}:
On peut borner $\sigma$ par:
$$
\frac{log(2)}{27} \leq \sigma \leq \frac{\log (18)}{6}.
$$
\end{theoreme}

\textbf{preuve:} D'apr\`{e}s la preuve de le th\'eor\`{e}me
(\ref{borner de card(m,r,o)}) on a le r\'{e}sultat suivant:
$$
\log (\mathcal{C}ard (M^n_2,r,o))\leq n^3 \log (18) / 6 .
$$
Alors
$$
\log (\mathcal{C}ard (M^n_2,r,o))/n^3\leq\log (18) / 6 .
$$
Le lemme (\ref{card(M,r,o)/n!,card(M,r),card(M,r,o)}) donne:
\begin{eqnarray*}
\log (\mathcal{C}ard (M^n_2,r))/n^3&\leq&\log (\mathcal{C}ard (M^n_2,r,o))/n^3\\
&\leq&\log (18) /6.
\end{eqnarray*}
Donc $\sigma\leq\log (18) / 6 .$\\
D'autre part, on a $2^{[n/3]^3}/n! \leq Card(M^n_2,r)$ voir le
lemme (\ref{borner de card(m,r,o)}) alors:
\begin{eqnarray*}
2^{[n/3]^3}/n! &\leq & Card(M^n_2,r)\\
\Rightarrow\qquad\qquad\,\,\,\, n^3/27\log(2)-\log(n!)&\leq
&\log(Card(M^n_2,r))\\
\Rightarrow\qquad\,\,\,\,\,\,\,\,\,\,\,\,\,\,
\log(2)/27-\log(n!)/n^3&\leq
&\log(Card(M^n_2,r))/n^3\\
\Rightarrow\qquad\,\, \log(2)/27-\lim_{n  \to
\infty}\log(n!)/n^3&\leq
&\sigma\\
\Rightarrow\qquad\qquad\qquad\qquad\,\, \log(2)/27-0&\leq &\sigma\\
\Rightarrow\qquad\qquad\qquad \qquad\qquad log(2)/27&\leq &\sigma.
\end{eqnarray*}

Donc,
$$
\frac{log(2)}{27} \leq \sigma \leq \frac{\log (18)}{6}.
$$

\end{document}